\documentclass{ert-l}

\newcommand{\kl}{l}
\newcommand{\Q}{\mathbb{Q}}
\newcommand{\C}{\mathcal{C}}
\newcommand{\R}{\mathbb{R}}
\newcommand{\Z}{\mathbb{Z}}
\newcommand{\N}{\mathbb{N}}
\newcommand{\x}{\mathfrak{x}}\newcommand{\y}{\mathfrak{y}}\newcommand{\z}{\mathfrak{z}}
\newcommand{\al}{\alpha}
\newcommand{\F}{\mathbb{F}_p}
\newcommand{\Fa}{\mathcal{F}}

\newcommand{\ol}{\mathfrak{o}}

\newcommand{\db}{{\bf b}}

\newcommand{\Dh}{D(H,K)} \newcommand{\Dg}{D(G,K)}\newcommand{\Dgr}{D_r(G,K)}
\newcommand{\Do}{D(H_0,K)}\newcommand{\Dgo}{D(G_0,K)} \newcommand{\Dgor}{D_r(G_0,K)}
\newcommand{\Dhr}{D_r(H,K)}\newcommand{\Dhrk}{D_{\bar{r}}}
\newcommand{\Dor}{D_r(H_0,K)}

\newcommand{\gr}{gr\dot{}_r\,}\newcommand{\gor}{gr\dot{}\,}\newcommand{\grd}{gr\dot{}_s\,}
\newcommand{\Co}{C^{an}(H_0,K)}\newcommand{\Cgo}{C^{an}(G_0,K)}

\newcommand{\Cg}{C^{an}(G,K)}

\newcommand{\B}{{\bf B}}

\newcommand{\Dgm}{D(G^m,K)}
\newcommand{\Dgom}{D(G^m_0,K)}

\newcommand{\lie}{\mathfrak{g}}

\newcommand{\nr}{\|.\|_{\bar{r}}}
\newcommand{\radius}{p^{-\frac{1}{p-1}}}
\newcommand{\stanh}{(\mathfrak{m} ^h)^d}

\newcommand{\Rm}{R^{(m)}}
\newcommand{\Rein}{R^{(1)}}
\newcommand{\Rmin}{R^{(m-1)}}
\newcommand{\bfh}{{\bf h}}
\newcommand{\Rep}{\mathcal{R}}
\newcommand{\T}{\Lambda}
\newcommand{\kap}{\kappa}

\newtheorem{theorem}{Theorem}[section]
\newtheorem{lemma}[theorem]{Lemma}
\newtheorem{proposition}[theorem]{Proposition}
\newtheorem{corollary}[theorem]{Corollary}
\theoremstyle{definition}

\theoremstyle{remark}

\numberwithin{equation}{section}

\begin{document}

\title[Auslander Regularity]{Auslander Regularity of $p$-adic Distribution Algebras}



\author[Tobias Schmidt]{Tobias Schmidt}
\address{Mathematisches Institut\\ Westf\"alische Wilhelms-Universit\"at
M\"unster\\ Einsteinstr. 62\\ D-48149 M\"unster, Germany}
\curraddr{D\'epartement de Math\'ematiques\\B\^atiment 425\\
Universit\'e Paris-Sud 11\\F-91405 Orsay Cedex, France}
\email{toschmid@math.uni-muenster.de}
\thanks{The author was supported by a grant within the DFG Graduiertenkolleg "Analytic Topology and Meta
Geometry" at M\"unster}

\subjclass[2000]{Primary  22E50; Secondary 11S99}

\date{July 3, 2007}

\dedicatory{}

\begin{abstract}
Given a compact $p$-adic Lie group over an arbitrary base field we
prove that its distribution algebra is Fr\'echet-Stein with
Auslander regular Banach algebras whose global dimensions are
bounded above by the dimension of the group. As an application, we
show that nonzero coadmissible modules coming from smooth or, more
general,~$U(\lie)$-finite representations have a maximal grade
number (codimension) equal to the dimension of the group.
\end{abstract}

\maketitle

\section{Introduction}
Given a locally $L$-analytic group $G$ where
$L\subseteq\mathbb{C}_p$ is a finite extension of $\Q_p$ \,P.
Schneider and J. Teitelbaum recently developed a systematic
framework to study locally analytic $G$-representations in
topological $p$-adic vector spaces (cf. [ST1-6]). At the center of
this theory lies a certain subcategory $\C_G$ of modules over the
$K$-valued locally analytic distribution algebra $\Dg$,
$K\subseteq\mathbb{C}_p$ being a complete and discretely valued
extension of $L$. The category $\C_G$ (the {\it coadmissible}
modules) is contravariantly equivalent to the category of
admissible locally analytic $G$-representations via the functor
"passage to the strong dual". The construction of $\C_G$ relies on
the fact that, when $G$ is compact, the algebra $\Dg$ is
Fr\'echet-Stein. The latter means that $\Dg$ equals a projective
limit of certain noetherian Banach algebras $\Dgr$ with flat
transition maps. Furthermore, Schneider and Teitelbaum prove
(\cite{ST5}, Thm. 8.9) that, in case $L=\Q_p$ and $G$ is compact,
the ring $\Dg$ is "almost" Auslander regular: the rings $\Dgr$ are
Auslander regular with a global dimension bounded above by the
dimension of the manifold $G$. This result allows to establish a
well-behaved dimension theory on $\C_G$ where the usual grade
number serves as a codimension function. On the other hand, it
allows to deduce important properties of the locally analytic
duality functor (in the sense of \cite{ST6}) e.g. its
involutivity.

The present work establishes the regularity property of $\Dg$ over
arbitrary base fields. Its main result is the

~\\ {\bf Theorem.}{\it~ Let $G$ be a compact locally $L$-analytic
group. Then $\Dg$ has the structure of a $K$-Fr\'echet-Stein
algebra where the corresponding Banach algebras are Auslander
regular rings whose global dimensions are bounded above by the
dimension of $G$.}

 ~\\We remark straightaway that the desired Banach algebras arise
in the same fashion as in \cite{ST5} i.e. they come as quotient
Banach algebras via the map $\Dgo\rightarrow\Dg$. Here, $G_0$
denotes the underlying locally $\Q_p$-analytic group and the map
is dual to embedding locally $L$-analytic functions into locally
$\Q_p$-analytic functions on $G$.

The brief outline of the paper is as follows. After reviewing the
notions of uniform pro-$p$ groups and Fr\'echet-Stein algebras in
our setting (sections 2-3) we begin by investigating certain
standard subgroups $H$ of $G$ serving as locally $L$-analytic
analogues of uniform groups (section 4). In section 5 we
generalize the filtration methods developed in \cite{ST5} to the
base field $L$. As a result we may deduce the regularity
properties for the algebras $\Dhr$ when $r$ is "sufficiently
small" and $H$ is as above. We then use the existence of a
well-behaved filtration of $\Dh$ by Fr\'echet-Stein subalgebras
(arising functorially from the lower $p$-series of $H$) to deduce
the main result (sections 6-7). We finish by describing the
resulting dimension theory on $\mathcal{C}_G$. As an application
we show that coadmissible modules coming from smooth or, more
general, $U(\lie)$-finite $G$-representations (as studied in
\cite{ST1}) are zero-dimensional.

~\\ {\it Notations.}~Throughout, $\Q_p\subseteq L\subseteq
K\subseteq \mathbb{C}_p$ is a chain of complete intermediate
fields where $L/\Q_p$ is finite and $K$ is discretely valued. If
not otherwise stated $G$ denotes a compact locally $L$-analytic
group, $\lie_L$ its Lie algebra and $\exp$ an exponential map. Let
$\kap:=1$ resp. $\kap:=2$ if $p$ is odd resp. even.

\section{Uniform pro-$p$ groups and $p$-valuations}

We begin by clarifying the relation between uniform pro-$p$ groups
(as introduced in \cite{DDMS}) and certain $p$-valued groups
introduced in \cite{ST5}. Let us first recall the definitions of
these groups and their basic properties.

Let $G$ be a topologically finitely generated pro-$p$ group.
Denoting by $G^l$ the subgroup of $G$ generated by $l$-th powers
the {\it lower $p$-series} $(P_i(G))_{i\geq 1}$ is inductively
defined via
\[P_1(G):=G,~P_{i+1}(G):=P_i(G)^p[P_i(G),G]\]
for all $i\geq 1$. The subgroups $P_i(G)$ are (topologically)
characteristic in $G$ and constitute a fundamental system of open
neighbourhoods for $1\in G$. The group $G$ is called {\it
powerful} if the commutator of $G$ is contained in $G^p$ resp.
$G^4$ (in case $p$ odd resp. even). In this case
$P_{i+1}(G)=P_i(G)^p=G^{p^{i}}$. Finally, $G$ is called {\it
uniform} if it is powerful and the index
$(P_{i}(G):P_{i+1}(G)),~i\geq 1$ does not depend on $i$. A uniform
group $G$ has a unique locally $\Q_p$-analytic structure: for any
minimal (ordered) set of topological generators $h_1,...,h_d$ the
map
\[(x_1,...,x_d)\mapsto h_1^{x_1}\cdot\cdot\cdot h_d^{x_d}\] is a bijective global chart
$\Z_p^d\rightarrow G$. In this situation, the subgroup $P_i(G)$
equals the image of $p^{i-1}\Z_p^d$ and is a uniform group itself.

~\\On the other hand, there is the notion of a $p$-{\it valuation}
$\omega$ on an abstract group $G$ introduced in \cite{Laz}. This
is a real valued function
\[\omega: G\setminus\{1\}\longrightarrow (1/(p-1),\infty)\] satisfying
\begin{enumerate}
  \item $\omega(gh^{-1}) \geq {\rm min~}(\omega(g),\omega(h)),$
  \item
 $ \omega(g^{-1}h^{-1}gh) \geq \omega(g)+\omega(h),$
  \item
 $ \omega(g^p) =\omega(g)+1$ \\
\end{enumerate}
for all $g,h\in G$. As usual one puts $\omega(1)=\infty$. A
$p$-valuation gives rise to a natural filtration of $G$ by
subgroups defining a topology on $G$. A $p$-valued group
$(G,\omega)$, complete with respect to this topology, is called
$p$-{\it saturated} if any $g\in G$ such that $\omega(g)>p/(p-1)$
is a $p$-th power.

Now let $G$ be a compact locally $\Q_p$-analytic group endowed
with a $p$-valuation $\omega$. It follows from \cite{Laz},
III.3.1.3/9 and III.3.2.1 that the topology on $G$ is defined by
$\omega$ whence $G$ is complete. Furthermore, $G$ admits an {\it
ordered basis}. This is an ordered set of topological generators
$h_1,...,h_d$ of $G$ such that the map $(x_1,...,x_d)\mapsto
h_1^{x_1}\cdot\cdot\cdot h_d^{x_d}$ is a bijective global chart
$\Z_p^d\rightarrow G$ and satisfies
\begin{equation}\label{sicher}\omega(h_1^{x_1}\cdot\cdot\cdot h_d^{x_d})=\min_{i=1,...,d}
(\omega(h_i)+v_p(x_i)).\end{equation} Here, $v_p$ is the $p$-adic
valuation on $\Z_p$. In \cite {ST5}, Sect.~4 the authors introduce
the class of compact locally $\Q_p$-analytic groups $G$ carrying a
$p$-valuation $\omega$ with ordered basis $h_1,...,h_d$ that
satisfy the following additional axiom

~\\
$\begin{array}{ll} {\rm (HYP)}
 & (G,\omega) {\rm ~is}~p-{\rm saturated~ and~ the~ ordered~ basis}~ h_1,...,h_d {\rm ~of~} $G$  \\
   & {\rm satisfies}~\omega(h_i)+\omega(h_j)>p/(p-1) {\rm ~for~ any~} 1\leq i\neq j\leq d. \\
\end{array}$
~\\~\\If $d=1$ the second condition is redundant.

We show that this class and the class of uniform groups are
closely related (comp. also [loc.cit.], remark after Lem. 4.3).
\begin{proposition}\label{pval}
Let $G$ be a uniform pro-$p$ group of dimension $d$. Then $G$ has
a $p$-valuation $\omega$ satisfying {\rm (HYP)}. It is given as
follows: for $g\in P_i(G)\setminus P_{i+1}(G)$ put $\omega(g):=i$
resp. $\omega(g):=i+1$ in case $p\neq 2$ resp. $p=2$. In
particular, $\omega$ is integrally valued. Any ordered system of
topological generators $h_1,...,h_d$ for $G$ is an ordered basis
for $\omega$ in the sense of {\rm (HYP)} with the property
\[\omega(h_1)=...=\omega(h_d)=\kap.\]
Conversely, if $p\neq 2$ then any compact locally $\Q_p$-analytic
group with a $p$-valuation satisfying {\rm (HYP)} is a uniform
pro-$p$ group.
\end{proposition}
\begin{proof}
We refer to \cite{DDMS} for all basic properties of uniform groups
that we use. Let $G$ be a uniform pro-$p$ group. Define $\omega$
as in the proposition. If $d>1$ then clearly
$\omega(h_i)+\omega(h_j)> p/(p-1)$ for $i\neq j$. Moreover,
$\omega(g)>p/(p-1)$ for $g\in G$ implies $g\in P_2(G)$ and the
group $P_2(G)$ consists (as a set) precisely of the $p$-th powers
of $G$. So for the first statement it remains to see that $\omega$
really is a $p$-valuation: the axiom 1. is clear. The map
$x\mapsto x^p$ is a bijection $P_i(G)/P_{i+1}(G)\rightarrow
P_{i+1}(G)/P_{i+2}(G)$ for all $i$ whence 3. Now the lower
$p$-series of any pro-$p$ group satisfies
\[[P_i(G),P_j(G)]\leq P_{i+j+\kap-1}(G)\] for all $i,j\geq 1$ which gives 2.
Finally, any ordered set of (topological) generators $h_1,...,h_d$
is an ordered basis and must lie in $P_1(G)\setminus P_2(G)$. In
particular, $\omega(h_i)=\kap$ for all $i$.

Conversely, assume $p\neq 2$ and let a compact locally
$\Q_p$-analytic group $G$ be given together with a $p$-valuation
$\omega$ satisfying (HYP). One has $\omega([g,h])>p/(p-1) $ for
all $g,h\in G$ according to {\rm (HYP)} and the equation
(\ref{sicher}) above. Since $G$ is $p$-saturated this implies
$[G,G]\subseteq G^p$ i.e. $G$ is powerful. It is also torsionfree
as there is no $p$-torsion according to 3. Hence, $G$ is uniform.
\end{proof}

\section{$\Dg$~as a Fr\'echet-Stein algebra}

We review the construction of a Fr\'echet-Stein structure on the
algebra $\Dg$ since this is central to our work. We refer to
\cite{ST5} for all details.

~\\Recall that a (two-sided) $K$-Fr\'echet algebra is called {\it
Fr\'echet-Stein} if there is a sequence $q_1\leq q_2 \leq...$ of
algebra seminorms on $A$ defining its Fr\'echet topology and such
that for all $n\in\N$ the completion $A_n$ of $A$ with respect to
$q_n$ is a noetherian $K$-Banach algebra and a flat
$A_{n+1}$-module via the natural map $A_{n+1}\rightarrow A_{n}$.
This applies to the algebra $\Dg$ of $K$-valued locally analytic
distributions on $G$. Recall that $\Dg$ equals the strong dual of
the locally convex vector space $\Cg$ of $K$-valued locally
$L$-analytic functions on $G$ equipped with the convolution
product. Let $G_0$ be the underlying locally $\Q_p$-analytic
group. Choose a normal open subgroup $H_0\subseteq G_0$ which is a
uniform pro-$p$ group. According to Prop. \ref{pval} the lower
$p$-series of $H_0$ induces an integrally valued $p$-valuation
$\omega$ on $H_0$ satisfying (HYP) with $h_1,...,h_d$ chosen to be
any (ordered) minimal set of generators of $H_0$. Thus, the
construction of \cite{ST5}, Sect. 4 for distribution algebras of
such $p$-valued groups applies: the bijective global chart
$\Z_p^d\rightarrow H_0$ for the manifold $H_0$ given by
\[(x_1,...,x_d)\mapsto h_1^{x_1}...h_d^{x_d}\] induces
a topological isomorphism $\Co\simeq C^{an}(\Z_p^d,K)$ on
$K$-valued locally analytic functions. In this isomorphism the
right-hand side is a space of classical Mahler series and the dual
isomorphism $D(H_0,K)\simeq D(\Z_p^d,K)$ therefore realizes $\Do$
as a space of noncommutative power series. More precisely, putting
$b_i:=h_i-1\in\Z[G],~\db^\al:=b_1^{\al_1}...b_d^{\al_d}$ for
$\al\in\N_0^d$ the Fr\'echet space $\Do$ equals all convergent
series
\[\lambda=\sum_{\al\in\N_0^d}d_\al\db^\al\]
with $d_\al\in K$ such that the set $\{|d_\al|r^{\kap
|\al|}\}_\al$ is bounded for all $0<r<1$. Here $|.|$ denotes the
normalized absolute value on $\mathbb{C}_p$. The value
$\lambda(f)\in K$ of such a series on a function $f\in\Co$ with
Mahler expansion
\[f({\bf x})=\sum_{\al\in\N_0^d}c_\al {{\bf x}\choose
{\al}},~c_\al\in K\] is given by $\sum_{\al\in\N_0^d}c_\al d_\al$.
The family of norms $\|.\|_r,~0<r<1$ defined via
\[\|\lambda\|_r:=\sup_\al |d_\al|r^{\kap|\al|}\]
defines the Fr\'echet topology on $\Do$. Restricting to the
subfamily $p^{-1}<r<1,~r\in p^\Q$ these norms are multiplicative
and the norm completions $\Dor$ are $K$-Banach algebras realizing
a Fr\'echet-Stein structure on $\Do$.

Choose representatives $g_1=1,...,g_r$ for the cosets in $G/H$ and
define on $\Dgo=\oplus_i\, \Do\,g_i$ the norms
$||\sum_i\lambda_i\,g_i||_r:=\max_i\,||\lambda_i||_r$. The
completions $\Dgor$ give suitable Banach algebras for $\Dgo$.

Finally, $\Dg$ is equipped with the quotient norms coming from the
map $\Dgo\rightarrow\Dg$. The latter arises as the transpose to
the embedding $\Cg\subseteq\Cgo$. Again, norm completions give
suitable Banach algebras.

We conclude with another important feature of $\Do$ in case of a
locally $\Q_p$-analytic group $H_0$ which is uniform. Each algebra
$\Dor$ carries the (separated and exhaustive) norm filtration
defined by the additive subgroups
\[\begin{array}{rl}
  F^s_r\Dor:= & \{\lambda\in\Dor,~||\lambda||_r\leq p^{-s}\}, \\
   &  \\
  F^{s+}_r\Dor:= & \{\lambda\in\Dor,~||\lambda||_r < p^{-s}\}
  \\
\end{array}
\] for $s\in\R$ with graded ring
\[\gr\Dor:=\oplus_{s\in\R}~gr_r^s\Dor\]
where $gr_r^s\Dor:=F^s_r\Dor/F^{s+}_r\Dor$. For a nonzero
$\lambda\in\Dor$ denote by deg$(\lambda)\in\R$ the {\it degree} of
$\lambda$ in the filtration. The {\it principal symbol}
$\sigma(\lambda)\neq 0$ of $\lambda$ is then given by \[\lambda +
F_r^{s+}\Dor\in gr_r^s\Dor\subseteq\gr\Dor\] where
$s=$deg$(\lambda)$. Note that $\gor K\simeq
k[\epsilon_0,\epsilon_0^{-1}]$ where $k$ denotes the residue field
and $\epsilon_0$ is the principal symbol of a prime element for
$K$.
\begin{theorem}\label{poly}
Mapping $\sigma(b_i)\mapsto X_i$ yields a $\gor K$-algebra
isomorphism of $\gr\Dor$ onto the polynomial ring $(\gor
K)[X_1,...,X_d]$.
\end{theorem}
\begin{proof} \cite{ST5}, Thm. 4.5.\end{proof}
We will prove an analogue over $L$ of this result in section 5.

\section{Results on certain standard groups}

Over the ground field $\Q_p$ the Fr\'echet-Stein theory of the
distribution algebra depends heavily on the presence of uniform
subgroups: every compact locally $\Q_p$-analytic group contains an
open normal uniform subgroup (\cite{DDMS}, Cor. 8.34). It is the
principal aim of this section to generalize this latter fact in a
suitable way to the ground field $L$.

~\\
We begin by reviewing the notion of a standard group. Let
$\mathfrak{m}$ be the maximal ideal in the valuation ring $\ol$ of
$L$. A locally $L$-analytic group is called {\it standard} of {\it
level} $h,~h\in\N$ if it admits a global chart onto
$\stanh\subseteq L^d$ such that the group operation is given by a
single power series without constant term and with coefficients in
$\ol^d$. A standard group of level $1$ will simply be called
standard. Any compact locally $L$-analytic group contains an open
subgroup which is standard (\cite{B-L}, III. 7.3 Thm. 4). In case
$L=\Q_p$ a locally $\Q_p$-analytic group is called {\it standard*}
if it admits a global chart onto $p^\kap\Z_p^d$ such that the
group operation is given by a single power series without constant
term and with coefficients in $\Z_p^d$. A standard* group is
uniform according to \cite{DDMS}, Thm.~8.31. We deduce a lemma
from the proof of this result.
\begin{lemma}\label{topgun} Suppose $G$ is a locally
$\Q_p$-analytic group which is standard* with respect to the
global chart $\psi: G\rightarrow p^\kap\Z_p^d$. Denote by $e_i\in
\Z_p^d$ the $i$-th unit vector. The elements
$h_i:=\psi^{-1}(p^\kap e_i)$ constitute a minimal set
$h_1,...,h_d$ of topological generators of $G$.
\end{lemma}
\begin{proof}
By the proof of [loc.cit.] $\psi$ induces a group isomorphism
\[G/P_2(G)\rightarrow p^\kap\Z_p^d/p^{\kap+1}\Z_p^d\] where $P_2(G)$
equals the Frattini subgroup of $G$. Hence, $h_1,...,h_d$ is a
minimal set of topological generators.
\end{proof}

Let $n=[L:\Q_p]$,~$e'$ be the ramification index of $L/\Q_p$ and
$u$ a uniformizer for $\ol$. Let $l\in\N,~l\geq 2$ and let $G$ be
locally $L$-analytic of dimension $d$.
\begin{lemma}\label{stanuni}
Suppose $G$ is standard of level $\kl e'$ with respect to a global
chart $\psi$. Then it is standard of level $1$ with respect to
$u^{1-\kl e'}\cdot\psi$. Its scalar restriction $G_0$ is standard*
with respect to $G_0\rightarrow (p^{\kap}\ol)^d\rightarrow
p^\kap\Z_p^{nd}$ where the first map equals
$p^{\kap-\kl}\cdot\psi$ and the second is induced by an arbitrary
choice of $\Z_p$-basis for $\ol$.
\end{lemma}
\begin{proof}
Let $\y_1,...,\y_d$ be the $L$-basis of $\lie_L$ induced by $\psi$
and regard $\psi$ as a map $G\rightarrow\Gamma$ where
$\Gamma=\oplus_j\mathfrak{m}^{\kl e'}\y_j\subseteq\lie_L$. By
standardness we have for $g,h\in G$ and
$\psi(g)=\sum_j\lambda_j\y_j, ~\psi(h)=\sum_j\mu_j\y_j,$ with
$\lambda_j,\mu_j\in\mathfrak{m}^{\kl e'}$ that
\[\psi(gh)=\sum_j
F_j(\lambda_1,...,\lambda_d,\mu_1,...,\mu_d)\y_j\] where
$F_j(X_1,...,X_d,Y_1,...,Y_d)\in\ol[[(X_s),(Y_s)]]$ without
constant term.

Let us first prove the statement concerning $G_0$: take as an
$L$-basis of $\lie_L$ the elements $\x'_j:=p^{\kl-\kap}\y_j.$ Then
$\Gamma=\oplus_j\mathfrak{m}^{\kl
e'}\y_j=\oplus_j\mathfrak{m}^{\kap e'}\x'_j$ and for $g,h\in G$
and $\psi(g)=\sum_i\lambda_j\x'_j, ~\psi(h)=\sum_j\mu_j\x'_j,$
with $\lambda_j,\mu_j\in\mathfrak{m}^{\kap e'}$ we have
\[\psi(gh)=\sum_j
F_j(\lambda_1p^{\kl-\kap},...,\lambda_dp^{\kl-\kap},\mu_1p^{\kl-\kap},...,\mu_dp^{\kl-\kap})\y_j.\]
Since $F_j\in\ol[[(X_s),(Y_s)]]$ has no constant term and $\kl\geq
2$ (i.e. $p^{\kl-\kap}\in\Z_p$) we get
\begin{equation}\label{becker}
\psi(gh)=\sum_j
F''_j(\lambda_1,...,\lambda_d,\mu_1,...,\mu_d)\x'_j\end{equation}
where $F''_j\in\ol[[(X_s),(Y_s)]]$ without constant term. By
definition $G$ is standard of level $\kap e'$ with respect to
\begin{equation}\label{chart}
G\stackrel{\psi}{\longrightarrow}
\Gamma=\oplus_j\mathfrak{m}^{\kap e'}\x'_j\longrightarrow
(\mathfrak{m}^{\kap e'})^d.\end{equation}

Now choosing a $\Z_p$-basis $v_1,...,v_n$ of $\ol$ yields
$\mathfrak{m}^{\kap e'}=\oplus_i\,p^\kap\Z_pv_i$. From
(\ref{chart}) together with the $\Q_p$-basis $\{v_i\x_j'\}$ of
$\lie_L$ we obtain the global chart
\[
G_0\stackrel{\psi}{\longrightarrow} \Gamma=\oplus_j\oplus_i
p^\kap\Z_pv_i\x'_j \longrightarrow (p^\kap\Z_p)^{nd}\] for $G_0$
which is as claimed. By (\ref{becker}) it follows for $g,h\in G_0$
and $\psi(g)=\sum_{ij}\lambda_{ij}v_i\x'_j,
~\psi(h)=\sum_{ij}\mu_{ij}v_i\x'_j,$ with
$\lambda_{ij},\mu_{ij}\in p^\kap\Z_p$ that
\[\begin{array}{rl}
   \psi(gh)= &\sum_j
F''_j((\sum_r\lambda_{rs}v_r)_s,(\sum_r\mu_{rs}v_r)_s)\x'_j \\
  &  \\
 =&  \sum_j
\sum_i G_{ij}((\lambda_{rs}),(\mu_{rs}))v_i\x'_j. \\
\end{array}\]
The $nd$ functions $G_{ij}$ are given by power series with
coefficients in $\Z_p$ and no constant term. Hence, $G_0$ is
standard* with a global chart as claimed.

It remains to see that $G$ is standard of level $1$ with respect
to the claimed chart. But this follows as above by taking
$\x_j:=u^{\kl e'-1}\y_j$ as $L$-basis of $\lie_L$.
\end{proof}
\begin{proposition}\label{normalstan}
Any compact locally $L$-analytic group $G$ contains an open normal
standard subgroup $H$ such that $H_0$ is standard*. The
corresponding global chart is given by
\[H\stackrel{\exp^{-1}}{\longrightarrow}\oplus_j\mathfrak{m}\x_j\longrightarrow \mathfrak{m}^d\]
where $\exp$ is the exponential map and $\x_1,...,\x_d$ a suitable
$L$-basis of $\lie_L$.
\end{proposition}
\begin{proof}
Choose a basis $\y'_1,...,\y'_d$ of $\lie_L$, endow $\lie_L$ with
the maximum norm and the space of $L$-linear maps on $\lie_L$ with
the operator norm. Scaling the basis we may assume that the
Hausdorff series converges on
\[\Lambda':=\oplus_j\mathfrak{m}\y'_j\] (viewed as a subset of $L^d$ via
$\y_1',...,\y_d'$) turning it into a locally $L$-analytic group.
We obtain an isomorphism of locally $L$-analytic groups $\exp:
\Lambda'\rightarrow G'$ onto an open subgroup $G'\subseteq G$.

Scaling the basis we may also assume that $|{\rm ad}~\x|<\radius$
for all $\x\in\Lambda'$ and that ${\rm Ad}(g)=\sum_{k\in\N_0}
1/k!~({\rm ad\,\x})^k$ for all $g=\exp(\x)\in G'$. Hence, $|{\rm
Ad}(g)|=1$ and $\Lambda'$ is Ad$(g)$-stable for all $g\in G'$.
Furthermore, let $\mathcal{R}'$ be a (finite) system of
representatives for the cosets $G/G'$ and put
\[\Lambda:=\cap_{g\in\mathcal{R}'}{\rm Ad}(g)\Lambda'.\]
Then $\Lambda$ is Ad$(g)$-stable for all $g\in G$. Now
$p^t\Lambda$ is a subgroup of $\Lambda'$ when $t\in\N$ is big
enough. To see this choose $t\in\N$ big enough such that
$g\exp(\x)g^{-1}=\exp({\rm Ad}(g)\x)$ holds for all $\x\in
p^t\Lambda'$ and for all $g\in\Rep'$. Then $p^t\Lambda$ is stable
under $*$, the group operation of $\Lambda'$. Indeed, if $\x,\y\in
p^t\Lambda$ and $g\in\Rep'$ then writing
$\x=$Ad$(g)\x',~\y=$Ad$(g)\y'$ with $\x',\y'\in p^t\Lambda'$ one
may compute \[\begin{array}{rl}
\exp(\x*\y)=\exp({\rm Ad}(g)\x')\exp({\rm Ad}(g)\y')=& g\exp(\x')\exp(\y')g^{-1} \\
  &  \\
   =&  \exp({\rm Ad}(g)(\x'*\y')). \\
\end{array}\]
Observe here, that we have $\x'*\y'\in p^t\Lambda'$ since
$p^t\Lambda'$ is a subgroup of $\Lambda'$. The calculation implies
that $\x*\y\in {\rm Ad}(g) p^t\Lambda'$ and since this holds for
any $g\in\Rep'$ we have $\x*\y\in p^t\Lambda$. Since
$p^t\Lambda=-p^t\Lambda$ and $0\in p^t\Lambda$ we obtain that
$p^t\Lambda$ really is a subgroup of $\Lambda'$ whence
\[M:=\exp(p^t\Lambda)\] is an open subgroup of $G$. Choose $\y_1,...,\y_d\in\lie_L$ such that
\[p^t\Lambda=\oplus_j\mathfrak{m}\y_j\] and consider the global
chart
\[\varphi: \exp(p^t\Lambda)\rightarrow
p^t\Lambda=\oplus_j\mathfrak{m}\y_j\rightarrow \mathfrak{m}^d.\]
We may assume that $M$ is standard with respect to $\varphi$
(otherwise pass to $\exp(\lambda^{-1}\mathfrak{m}^d)$ for
sufficiently big $\lambda\in L$ and use the chart
$\lambda\cdot\varphi$, see (the proof of) \cite{B-L}, III.7.3
Thm.~4). In this situation consider
\[\Gamma:=\oplus_j\mathfrak{m}^{\kl e'}\y_j\subseteq
\oplus_j\mathfrak{m}\y_j=p^t\Lambda\] for suitable
$\kl\in\N,~\kl\geq 2$ which will be chosen conveniently in the
following. In the rest of the proof we show that the open subgroup
\[H:=\exp(\Gamma)\] of $G$ will satisfy our requirements.

Let us first check normality: $H$ is normal in $M$ by \cite{B-L},
III.7.4 Prop. 6. Let $\Rep$ be a (finite) system of
representatives for the cosets in $G/M$. Enlarging $\kl$ if
necessary we may assume that $g\exp(\x)g^{-1}=\exp({\rm Ad}(g)\x)$
holds all $\x\in\Gamma,~g\in\Rep$. Now $\Lambda$ is Ad$(g)$-stable
for all $g\in G$ and hence, so is $\Gamma=u^{\kl e'-1}p^t\Lambda$
($u$ a uniformizer for $\ol$) since each Ad$(g)$ is $L$-linear.
Thus, $H=\exp(\Gamma)$ is stable under conjugation with elements
from $\Rep$ and therefore it is normal in $G$.

Now since $M$ is standard with respect to $\varphi$ it follows
from [loc.cit.] that $H$ is standard of level $\kl e'$ with
respect to $\varphi|_{_H}$. Applying Lem. \ref{stanuni} we see
that $H$ is also standard of level $1$ with a global chart given
by
\[H\stackrel{\exp^{-1}}{\longrightarrow}\Gamma=\oplus_j\mathfrak{m}\x_j\longrightarrow
\mathfrak{m}^d\] where $\x_j:=u^{\kl e'-1}\y_j$. Thus, the
$L$-basis $\x_1,...,\x_d$ is as desired. Finally, Lem.
\ref{stanuni} also implies that the restricted group $H_0$ is
standard*.
\end{proof}
Given a locally $L$-analytic group $G$ a choice of $L$-basis
$\x_1,...,\x_d$ of $\lie_L$ gives rise to the map
\[\theta_L: (\sum_jx_j\x_j)\mapsto
\exp(x_1\x_1)\cdot\cdot\cdot\exp(x_d\x_d)\in G\] defined on an
open subset of $\lie_L=\oplus_j L\x_j$ containing $0$. It is
locally $L$-analytic and \'etale at $0$ and is called a {\it
system of coordinates of the second kind} associated to the
decomposition $\lie_L=\oplus_j L\x_j$ (\cite{B-L}, III.4.3
Prop.~3).

Identifying the Lie algebras $\lie_L$ resp. $\lie_{\Q_p}$ of $G$
resp. $G_0$ over $\Q_p$ (according to \cite{B-VAR}, 5.14.5) we
consider the following condition on $G$:

~\\{\it Condition {\rm (L)}: There is an $L$-basis $\x_1,...,\x_d$
of $\lie_L$ and a $\Z_p$-basis $v_1,...,v_n$ of $\ol$ such that
the system of coordinates of the second kind $\theta_{\Q_p}$
induced by the decomposition $\lie_{\Q_p}=\oplus_j\oplus_i \Q_p
v_i\x_j$ gives an isomorphism of locally $\Q_p$-analytic manifolds
\[\theta_{\Q_p}:\oplus_j\oplus_i\Z_pv_i\x_j\longrightarrow G_0.\]
The exponential satisfies $\exp (\lambda\cdot v_i\x_j)=\exp
(v_i\x_j)^\lambda$ for all $\lambda\in\Z$.}

~\\Note that if $G$ is pro-$p$ and satisfies (L) with suitable
bases $v_i$ and $\x_j$ then $\exp(\lambda\cdot
v_i\x_j)=\exp(v_i\x_j)^\lambda$ for all $\lambda\in\Z$ extends to
$\Z_p$-powers and so the elements
\[h_{ij}:=\theta_{\Q_p}(v_i\x_j)=\exp(v_i \x_j)\] are a minimal ordered system
of topological generators for $G_0$.

For our purposes locally $L$-analytic groups that are uniform and
satisfy the above condition (such as the additive group $\ol$)
will play the role of uniform locally $\Q_p$-analytic groups. We
show that there are enough of them.
\begin{corollary}\label{axiomL} Any compact
locally $L$-analytic group $G$ has a fundamental system of open
normal subgroups $H$ such that $H_0$ is uniform and satisfies {\rm
(L)}.
\end{corollary}
\begin{proof}
It suffices to show that the group $H$ constructed in the proof of
the last proposition satisfies (L). This is because $H\subseteq
G'$ (in the notation of this proof) and by construction, the open
subgroup $G'$ of $G$ can be chosen as small as desired. We use the
notation of this proof.

According to it (and in connection with Lem. \ref{stanuni}) $H_0$
is standard* with respect to the bijective global chart
\[\psi: H_0 \stackrel{\exp^{-1}}{\longrightarrow}
 \Gamma=\oplus_j\oplus_ip^\kap\Z_pv_i\x_j'\longrightarrow p^\kap\Z_p^{nd}.\]
Here, $\x'_1,...,\x'_d$ and $v_1,...,v_n$ are bases of $\lie_L$
resp. $\ol$. Denoting by $e_{ij}$ the $ij$-th unit vector in
$\Z_p^{nd}$ put
\[h_{ij}:=\psi^{-1}(p^\kap e_{ij})=\exp (p^\kap v_i\x_j')\]
for $i=1,...,n,~j=1,...,d$. Lem. \ref{topgun} implies that these
$nd$ elements are a minimal system $h_{11},h_{21},...,h_{nd}$ of
topological generators for $H_0$. Putting $\z_j:=p^\kap\x'_j$ the
map
\[\theta_{\Q_p}:\oplus_j\oplus_i\Z_pv_i\z_j\longrightarrow
H_0,~\sum_{ij}\lambda_{ij}v_i\z_j\mapsto \prod_j\prod_i
h_{ij}^{\lambda_{ij}}\] is a locally $\Q_p$-analytic isomorphism
since $H_0$ is uniform (section 2). Because of
$h_{ij}=\exp(v_i\z_j)$ this map is by definition the system of
coordinates of the second kind induced by the decomposition
$\lie_{\Q_p}=\oplus_j\oplus_i\Q_pv_i\z_j$. The last condition of
{\rm (L)} follows from this as well.
\end{proof}
\begin{corollary}\label{axiomLm}
Suppose $H$ is an open normal subgroup of $G$ such that $H_0$ is
uniform satisfying {\rm (L)}. Then any member in the lower
$p$-series of $H_0$ is an open normal subgroup of $G$ whose
restriction is uniform satisfying {\rm (L)}.
\end{corollary}
\begin{proof}
Let $H^m=\{x^{p^m}, x\in H\},~m\in\N_0$ be a member in the lower
$p$-series of $H_0$ with locally $L$-analytic structure as an open
subgroup of $G$. It is thus open normal in $G$. Let
$\x_1,...,\x_d$ resp. $v_1,...,v_n$ be bases of $\lie_L$ resp.
$\ol$ such that the induced coordinate system $\theta_{\Q_p}$
realizes condition (L) for $H_0$. It immediately follows that
$H^m_0$ satisfies (L) with respect to the bases
$p^m\x_1,...,p^m\x_d$ and $v_1,...,v_n$.
\end{proof}

\section{Filtrations}

We wish to extend the filtration methods developed in \cite{ST5}
over $\Q_p$ to distribution algebras of certain locally
$L$-analytic groups. In particular, we aim at an analogue over $L$
of Thm. \ref{poly}. The most natural way to do this is to realize
$\Dg$ as a quotient of $\Dgo$ and study the induced filtrations.

Let $G$ be as usual a compact locally $L$-analytic group. Via the
identification $\lie_L\simeq\lie_{\Q_p}$ over $\Q_p$ we view
$\exp$ as an exponential map for $G_0$ as well. Recall that
$\lie_L$ acts on $C^{an}(G_0,K)$ via differential operators in the
usual way, i.e. \[\x f (g):=\frac{d}{dt}
\,f(\exp(-t\x)g)\,|_{t=0}.\] This gives an $L$-linear inclusion
$L\otimes_{\Q_p}\lie_L\subseteq\Dgo$ via associating to
$1\otimes\x$ the functional $f\mapsto (-\x)f(1)$. The kernel of
the quotient map $\Dgo\rightarrow \Dg$
\[I(G_0,K):=\{\lambda\in\Dgo: \lambda|_{_{C^{an}(G,K)}}\equiv 0\}\]
is a two-sided and closed ideal.
\begin{lemma}
As a right ideal $I(G_0,K)$ is finitely generated by the elements
\[F_{ij}:=1\otimes v_i\x_j-v_i\otimes\x_j\in\Dgo.\]
Here, $\x_1,...,\x_d$ is any $L$-basis of $\lie_L$ and
$v_1,...,v_n$ is any $\Q_p$-basis of $L$.
\end{lemma}
\begin{proof}
For sake of clarity denote for the moment by $\delta_g$ the image
of $g\in G$ in $\Dgo$. Since $K[G]\subseteq\Dgo$ is dense it
suffices to show that $I:=I(G_0,K)$ equals the closure of the
$K$-vector space generated by the products
\[F_{ij}\delta_g\in\Dgo\] where $i=1,...,n,~j=1,...,d,~g\in G$. Let $W$ denote the closure of the
described vector space. According \cite{Ko}, Lem. $1.3.2$ the
subspace $\Cg\subseteq\Cgo$ equals
\[\{f\in\Cgo: (1\otimes v_i\x_j-v_i\otimes\x_j)f=0 {\rm~in~}\Cgo
{\rm~for~all~}i,j\}\] and so $F_{ij}\delta_g\in I$ whence
$W\subseteq I$. Choose a continuous functional $\phi$ on $I$
vanishing on the subspace $W$. We show $\phi=0$ whence $W=I$ by
Hahn-Banach. For the strong dual $I'_b$ we have the canonical
isomorphism $I'_b\simeq C^{an}(G_0,K)/\Cg$ (\cite{B-TVS}, Cor.
IV.2.2) whence $\phi=\bar{f}$ for some $\bar{f}\in\Cgo/\Cg$. Then
$\phi$ vanishing on $W$ implies
\[
F_{ij}\delta_g (f) =0.
\]
Recalling that multiplication in $\Dgo$ is convolution one
computes
\[\begin{array}{ccl}
  0 & = &-F_{ij}\delta_g (f)
  \\&&\\
   & = &\delta_g(g'\longrightarrow -F_{ij}(g''\longrightarrow f(g''g')))
   \\&&\\
   & = & \delta_g((1\otimes v_i\x_j)f-(v_i\otimes\x_j) f)=(1\otimes v_i\x_j-v_i\otimes\x_j)f(g)\end{array}.\]
Since $i,j$ and $g\in G$ were arbitrary we have $f\in\Cg$ which
means $\phi=0$.
\end{proof}

Remark: If $I_r(G_0,K)$ denotes the norm closure inside a Banach
algebra $\Dgor$ we have $I_r(G_0,K)=I(G_0,K)\Dgor$ whence
$I_r(G_0,K)$ is also finitely generated (as a right ideal of
$\Dgor$) by the $F_{ij}$.

~\\From now on we additionally assume that $G$ is uniform and
satisfies (L). Let $\x_1,...,\x_d$ resp. $~v_1,...,v_n$ be
corresponding bases of $\lie_L$ resp. $\ol$. We may and will
arrange $v_1\in\Z_p$ and thus $F_{1j}=0$ for all $j$. Put
$h_{ij}:=\exp (v_i\x_j)$ for $i,j\geq 1$. The elements
$h_{11},h_{21},...,h_{nd}$ constitute a minimal ordered set of
topological generators for the uniform group $G_0$ and the global
chart
\[ G_0\stackrel{\theta_{\Q_p}^{-1}}{\longrightarrow}
\oplus_j\oplus_{i}\Z_pv_i\x_j\longrightarrow \Z_p^{nd}\] induces
expansions for the elements of $\Dgo$ in the monomials
$\db^\al=b_{11}^{\al_{11}}b_{21}^{\al_{21}}...b_{nd}^{\al_{nd}},~\al\in\N_0^{nd}$
where $b_{ij}:=h_{ij}-1\in\Z[G]\subseteq\Dgo$ (section 3).
Finally, let $\log(1+X):=\sum_{k\geq 1}(-1)^{k-1}X^k/k\in\Q[[X]]$.
\begin{lemma}\label{expansions}
We have $F_{ij}=\log(1+b_{ij})-v_i\log(1+b_{1j}).$
\end{lemma}
\begin{proof}
Let $f\in\Cgo$. It has a Mahler expansion of the form
\[(f\circ\theta_{\Q_p})(\sum_{ij}x_{ij}v_i\x_j)=\sum_{\al}c_\al {{\bf x}\choose \al}\]
for all ${\bf x}:=(x_{11},x_{21},...,x_{nd})\in\Z_p^{nd}$ with
coefficients $c_\al\in K$. Since $\theta_{\Q_p}$ restricted to the
direct summand $\Q_pv_i\x_j$ equals $\exp$ we may explicitly
calculate the values $(v_i\otimes\x_j)(f)$ and $(1\otimes
v_i\x_j)(f)$. Then a comparison of coefficients yields the claim.
\end{proof}
Now let $||.||_r$ denote a norm on $\Dgo$ and $\nr$ the quotient
norm on $\Dg$. Consider the exact sequence
\[0\longrightarrow I_r(G_0,K)\longrightarrow \Dgor\longrightarrow\Dgr
\longrightarrow 0\] of filtered $\Dgor$-modules. By our
assumptions ($r\in p^\Q$ and $K$ discretely valued) the
filtrations in question are essentially indexed by $\Z$ and so the
functor $\gr$ is exact whence
\[\gr\Dgr\simeq\gr\Dgor/\gr I_r(G_0,K)\]
canonically. Recall that $\gor\Dgor=(\gor K)[X_{11},...,X_{nd}]$
(Thm. \ref{poly}) and so we aim at calculating $\gr I_r(G_0,K)$.

Abbreviate $D_r:=D_r(G_0,K),~I_r:=I_r(G_0,K), ~D_{\bar{r}}:=\Dgr$
for the rest of this section.

We start by computing the principal symbols $\sigma(F_{ij})$ and
then deduce that they generate $\gr I_r$. Recall that a real
number $0<s<1$ is called {\it not critical} for the series
$\log(1+X)$ if the supremum
\[|\log(1+X)|_s:=\sup_{k\in\N} \,|\frac{1}{k}X^k|_s=\sup_{k\in\N} \,|\frac{1}{k}|\,s^k\]
is obtained at a single monomial of the series $\log(1+X)$. This
monomial is called {\it dominant}. We assume for the rest of this
section that $r^\kap$ is not critical for $\log (1+X)$.
\begin{lemma}\label{noc}\label{symbol}
We have
\[\sigma(F_{ij})=\epsilon^{-h}\,X_{ij}^{p^h}-\bar{v}_i\epsilon^{-h}\,X_{1j}^{p^h}\in\gr
D_r\] with $h\in\N_0$ depending only on $r^\kap$,~$\bar{v}_i\in k$
is the residue class of $v_i$ and $\epsilon:=\sigma(p)\in\gor K$.
For all $r^\kap<\radius$ one has $h=0$ i.e.
$\sigma(F_{ij})=X_{ij}-\bar{v}_iX_{1j}.$
\end{lemma}
\begin{proof}
Since $r^\kap$ is not critical there is a dominant monomial of
$\log(1+X)$ with respect to $|.|_{r^\kap}$ whose index depends
only on $r^\kap$ and is easily seen to be a $p$-power, say $p^h$.
Using Lem. \ref{expansions} the claim follows by direct
calculation.
\end{proof}
\begin{corollary}\label{klinear}
The $nd-d$ elements $F_{ij},~i\neq 1$ are orthogonal in $D_r$ i.e.
one has for arbitrary $c_{ij}\in K$ that
\[ ||\sum_{ij}\,c_{ij}F_{ij}||_r=\max_{ij}||c_{ij}F_{ij}||_r.\]
\end{corollary}
\begin{proof}
We may assume (via leaving away possible summands) that
$c_{ij}\neq 0$ for all $ij$ and that
$||c_{ij}F_{ij}||_r=:p^{-s},~s\in\R$ is a constant for all $ij$.
According to the above lemma the elements $\sigma(F_{ij})$
generate a free $\gor K$-submodule inside $\gr D_r$. We therefore
get
\[0\neq
\sum_{ij}\sigma(c_{ij})\sigma(F_{ij})=\sum_{ij}\,\sigma(c_{ij}F_{ij})=\sum_{ij}\,c_{ij}F_{ij}~~~{\rm
mod~}F_r^{s+}D_r\] and so $||\sum_{ij}\,c_{ij}F_{ij}||_r=p^{-s}$.
\end{proof}

For the rest of this section write $\Fa$ for the set of elements
$F_{ij},~i\neq 1$. Put $s_0:=$deg$(F)$, some $F\in\Fa$. By Lem.
\ref{expansions} the value $s_0$ is independent of the choice of
$F$.
\begin{proposition}\label{goody} \label{gensym}For any given $s\in\R$ one has
\[F^s_rI_r=\sum_{F\in\Fa}F~F^{s-s_0}_rD_r\]
as additive groups. In particular, the graded ideal $\gr I_r$ is
generated by the principal symbols $\sigma(F),~F\in\Fa$.
\end{proposition}
\begin{proof}
Write $\tilde{F}_r^sI_r:=\sum_{F\in\Fa}F~F^{s-s_0}_rD_r$ and
similarly for $s+$ instead of $s$. The other inclusion being
trivial we prove $F^s_rI_r\subseteq\tilde{F}_r^sI_r$ for all $s$.
This reduces to show
\begin{equation}\label{key} \tilde{F}_r^sI_r\cap
F_r^{s+}I_r\subseteq \tilde{F}_r^{s+}I_r.
\end{equation} Indeed, if
$\lambda\in F_r^sI_r$, let $s'\leq s$ be a number such that
$\lambda\in\tilde{F}^{s'}_rI_r$. Since $r\in p^\Q$ and $K$ is
discretely valued the filtration $F_r\dot{} I_r$ is essentially
indexed by $\Z$. Using (\ref{key}) we see that we may choose
$s'=s$.

So let us prove (\ref{key}). Let $\lambda\in\tilde{F}_r^sI_r\cap
F_r^{s+}I_r$. Of all representations $\lambda=\sum_{k}
F_k\lambda_k\in F^{s+}_rI_r$ with $F_k\in\Fa$, $F_k\lambda_k\in
F^s_rI_r$ take one with minimal number $t$ of nonzero summands
$F_k\lambda_k$. It suffices to show that if $t>1$ then, modulo a
term $T$ in $\tilde{F}_r^{s+}I_r$, there is a representation with
$t-1$ nonzero summands. Putting summands contained in
$F_r^{s+}I_r$ into $T$ we may assume that deg$(F_k\lambda_k)=s$
for all $k$. Then the hypothesis implies $\sum_k
\sigma(F_k)\sigma(\lambda_k)=0$. Using Lem. \ref{symbol} and
\cite{Ka}, 3.1, Ex. 12 (c) it easily follows that
$\sigma(\lambda_k)$ is contained in the ideal generated by all
$\sigma(F_{k'}),~k'\neq k$. We now assume $t=2$ (the general case
follows in exactly the same fashion but with more notation). The
element $\sigma(\lambda_1)$ thus lies in the ideal generated by
$\sigma(F_2)$ whence $\lambda_1=F_2\lambda'_2+R$ with $\lambda'_2,
R\in D_r$ and $||R||_r<||\lambda_1||_r$. Hence
\[\lambda=F_1(F_2\lambda_2'+R)+F_2\lambda_2.\]
Now $||F_1F_2-F_2F_1||_r<||F_1F_2||_r$ since $\gr D_r$ is
commutative. One the other hand, the $nd-d$ elements in $\Fa$
generate (over $L$) the kernel of the Lie algebra map
$L\otimes_{\Q_p}\lie_{\Q_p}\rightarrow\lie_L,~a\otimes\x\mapsto
a\x$. It follows that $F_1F_2-F_2F_1=\sum_{F\in\Fa} Fc_F$ with
$c_F\in L$ and so, together with Cor. \ref{klinear}, we obtain
$|c_F|<||F_1||_r$ for all $F\in\Fa$. All in all we have $
\lambda\equiv F_2(F_1\lambda_2'+\lambda_2)$ modulo terms in
$\tilde{F}_r^{s+}I_r$ and $F_2(F_1\lambda_2'+\lambda_2)\in
F^s_rI_r$.
\end{proof}
Summarizing the results obtained so far yields the
\begin{proposition}\label{quotient1}\label{quotient2}
There is an isomorphism of $\gor K$-algebras
\[\gr\Dhrk\stackrel{\sim}\longrightarrow(\gor
K)[X_{11},...,X_{nd}]/(\{X_{ij}^{p^h}-\bar{v}_iX_{1j}^{p^h}\}_{i\geq
2,j\geq 1})\] where $h\in\N_0$ depends only on $r^\kap$. If
$r^\kap<\radius$ then $h=0$ and
\[\gr\Dhrk\simeq(\gor K)[X_{11},...,X_{1d}]\] as $\gor K$-algebras where the isomorphism is obtained by
the first map composed with the algebra homomorphism induced by
$X_{ij}\mapsto \bar{v}_iX_{1j}$ for all $i\geq 2,j\geq 1$. In this
case the norm $\nr$ is multiplicative and $\Dhrk$ is an integral
domain.
\end{proposition}
Since $\Dg$ embeds in any completion $\Dhrk$ there is the
\begin{corollary}
The ring $\Dg$ is an integral domain.
\end{corollary}
This generalizes the case $G:=\ol$ appearing in \cite{ST2}, Cor.
3.7.

Recall that $D_r$ equals a noncommutative power series ring in the
elements $\db^\al,~\al\in\N_0^{nd}$ where $b_{ij}=h_{ij}-1$.
Viewing the elements $\db^\al$ in $\Dhrk$ via $K[G]\subseteq\Dhrk$
we may consider their principal symbols inside $\gr\Dhrk$. Tracing
through the definitions of the maps involved yields the
\begin{corollary}\label{isoconnection}
Let $r^\kap<\radius$. The isomorphism \[\gr\Dhrk\rightarrow(\gor
K)[X_{11},...,X_{1d}]\] maps $\sigma(b_{1j})\mapsto X_{1j}$.
\end{corollary}
For an index $\beta\in\N_0^d$ (consisting of $d$ components!) let
us agree that $\db^\beta$ denotes the element
$b_{11}^{\beta_1}...b_{1d}^{\beta_d}\in K[G]$.
\begin{proposition}\label{result1b}
Let $r^\kap<\radius$. Every $\lambda\in\Dhrk$ has a convergent
expansion
\[\lambda=\sum_{\beta\in\N_0^d}d_\beta\db^\beta\]
with uniquely determined $d_\beta\in K$. Furthermore,
\[||\lambda||_{\bar{r}}=\sup_\beta||d_\beta\db^\beta||_{\bar{r}}=\sup_\beta|d_\beta|r^{\kap\,|\beta|}.\]
\end{proposition}
\begin{proof}
The norm $\nr$ is multiplicative and so $\gr\Dhrk$ is freely
generated as $\gor K$-module by the symbols
$\{\sigma(\db^\beta)\}_{\beta\in\N_0^d}$ according to Cor.
\ref{isoconnection}. This implies that the monomials $\db^\beta$
generate a dense $K$-subspace in $\Dhrk$ such that
$||\lambda||_{\bar{r}}=\max_\beta ||d_\beta\db^\beta||_{\bar{r}}$
for any finite sum $\lambda=\sum_\beta d_\beta\db^\beta$ out of
this submodule. Finally, Cor. \ref{isoconnection} yields that
$\sigma(b_{1j})\in\gr D_r$ is mapped to
$\sigma(b_{1j})\in\gr\Dhrk$ by the canonical map $\gr
D_r\rightarrow\gr D_{\bar{r}}$ whence
$||b_{1j}||_r=||b_{1j}||_{\bar{r}}$.
\end{proof}
Remark: When $r^\kap>\radius$ it is easy to see that $\gr\Dhrk$
has nonzero nilpotents, in contrast to Cor. \ref{isoconnection}.

\section{Lower $p$-series subalgebras}

In this section we prepare the proof of the main result.

~\\To begin with let $G$ be an arbitrary uniform locally
$\Q_p$-analytic group of dimension dim$_{\Q_p}G=d$. Fix a sequence
of topological generators $h_1,...,h_d$. Let
\[G=G^0\supseteq G^1\supseteq G^2\supseteq...\] be its lower
$p$-series. By functoriality we obtain a series of subalgebras
\[\Dg=D(G^0,K)\supseteq D(G^1,K)\supseteq D(G^2,K)\supseteq...\] where
the inclusions are topological embeddings with respect to
Fr\'echet topologies (\cite{Ko}, Prop. 1.1.2). If $n\geq m$ then
$\Dgm$ is a finitely generated free (left or right) module over
$D(G^n,K)$ on a basis any set of coset representatives for
$G^m/G^n$. Since each $G^m$ is a uniform pro-$p$ group, generated
by $h_1^{p^m},...,h_d^{p^m}$, the algebra $\Dgm$ is
Fr\'echet-Stein and equals a noncommutative power series ring in
the elements $h_1^{p^m}-1,...,h_d^{p^m}-1$. We denote the induced
set of norms on $\Dgm$ by $||.||_r^{(m)},~ p^{-1}<1<1,~r\in p^\Q$.
Abbreviate $||.||_r:=||.||_r^{(0)}$.

For clarity we recall two simple definitions. If $I$ denotes a
countable index set a family of pairwise different nonzero
elements $(v_i)_{i\in I}$ in a normed $K$-vectorspace $(V, ||.||)$
is called {\it orthogonal} in  $V$ if one has $||v||=\max_I
|c_i|\,||v_i||$ for any convergent series $v=\sum_I\,c_iv_i,~
c_i\in K$. The family is called an {\it orthogonal basis} if any
element of $V$ can be written as such a convergent series.
\begin{lemma}\label{Frommer2}
Let $\Lambda=\{\lambda_i\}_{i\in I}$ be a countable family of
pairwise different nonzero elements of $\Dgr$. For any
 $\lambda_i\in\Lambda$ expand $\lambda_i=\sum_{\al\in\N_0^d} d_{i,\al}\db^\al$ and choose $\gamma_{i}\in\N_0^d$ with
 $||\lambda_i||_r=||d_{i,\gamma_{i}}\db^{\gamma_{i}}||_r$.
Suppose that for each $\lambda_i$ the index $\gamma_{i}$ is
uniquely determined. If the map $\iota:
\Lambda\rightarrow\N_0^d,~\lambda_i\mapsto \gamma_{i}$ is
injective resp. bijective the system $\T$ is orthogonal resp. an
orthogonal basis in $\Dgr$.
\end{lemma}
\begin{proof}
This follows easily from \cite{F}, Lem. 1.4.1/2.
\end{proof}
\begin{proposition}\label{restriction}
Fix $m\geq 1$ in $\N$ and suppose $r^{\kappa(p^m-1)}>p^{-1}$. Then
$\|.\|_r$ restricts to $\|.\|_{r'}^{(m)}$ on the subring
$\Dgm\subseteq\Dg$ where $r'=r^{p^m}$.
\end{proposition}
\begin{proof}
Abbreviate $R:=\Dg,~R^{(m)}:=\Dgm$ and use induction on $m$. Let
$m=1$. Put
$b'_i:=h_i^p-1,~\db'^\al:=b_1'^{\al_1}...b_d'^{\al_d},~\al\in\N_0^d$
and consider an arbitrary element
\[
\lambda=\sum_{\al\in\N_0^d}d_\al\db'^\al\]
 in $R^{(1)}$. To verify that $\|.\|_r$ restricts to $\|.\|^{(1)}_{r^p}$ on
$R^{(1)}$ amounts to show, by definition of the norms
$||.||_s^{(1)}$, that
$\|\lambda\|_r=\sup_\al|d_\al|r^{p\kap|\al|}.$ To do this we first
compute the norm $\|\db'^\al\|_r$. We have
\[b'_i=h_i^p-1=(b_i+1)^p-1=b_i^p+\sum_{k=1,...,p-1}{p\choose
k}b_i^k.\] By $r^{\kap(p-1)}>p^{-1}$ it follows
\[\|b_i^p\|_r=r^{\kap p}>p^{-1}r^\kap\geq p^{-1}r^{\kap k}=|{p\choose k}|r^{\kap k}=\|{p\choose
k}b_i^k\|_r\]for all $k=1,...,p-1$ and so $\|b'_i\|_r=r^{\kap p}$.
Hence, $\|\db'^\al\|_r=\prod_i\|b'_i\|_r^{\al_i}=r^{p\kap|\al|}.$

Since the Fr\'echet topology on $R^{(1)}$ is stronger than the
induced $||.||_r$-topology we may view $\lambda=\sum_\al
d_\al\db'^\al$ as a convergent series in the $K$-Banach space
$(\Dgr,||.||_r)$. Then we are reduced to show that the family of
pairwise different nonzero elements
$\Lambda:=\{\db'^\al\}_{\al\in\N_0^d}$ is orthogonal in
$(\Dgr,||.||_r)$. But this follows from Lem. \ref{Frommer2}: by
our calculations $\gamma_{\al}$ is uniquely determined and equals
$p\al$. Hence, $\iota: \db'^\al\mapsto p\al$ is injective.

Now let $m>1$ and assume that the result holds true
 for all numbers strictly smaller than $m$ and that
 $r^{\kap({p^m}-1)}>p^{-1}$. Since the latter implies
$r^{\kap(p^{m-1}-1)}>p^{-1}$ the induction
 hypothesis shows that $\|.\|_r$ on $R$ restricts to
 $\|.\|^{(m-1)}_{r^{p^{m-1}}}$ on $R^{(m-1)}$.

 Now $G^m$ appears also as first step in the lower $p$-series
 of $G^{m-1}$. Hence, the induction
 hypothesis applies to these two groups: for $0<s<1$, every $\|.\|^{(m-1)}_s$ on
 $R^{(m-1)}$ restricts to $\|.\|^{(m)}_{s^p}$ on $R^{(m)}$ as long as
 $s^{\kap p}>s^\kap p^{-1}$. We choose $s=r^{p^{m-1}}$. Then
 $s^{\kap p}=r^{\kap p^m}>r^\kap p^{-1}> s^\kap p^{-1}$ and so
$\|.\|^{(m-1)}_{r^{p^{m-1}}}$
 restricts to $\|.\|^{(m)}_{r^{p^m}}$ on $R^{(m)}$.
\end{proof}
Denote by $D_{(r)}(G^m,K)\subseteq\Dgr$ the norm closure of $\Dgm$
inside the Banach algebra $\Dgr$.
\begin{lemma}
Assume $r^{\kap(p^m-1)}>p^{-1}$. Then $\Dgr$ is a finite free
(left or right) $D_{(r)}(G^m,K)$-module on a basis any system of
coset representatives for $G/G^m$.
\end{lemma}
\begin{proof}
Since the argument for right modules is the same we only prove the
statement concerning left modules.

Put
$R^{(m)}:=\Dgm,~\Rm_r:=D_{(r)}(G^m,K),~R:=R^{(0)},~R_r:=R^{(0)}_r$.
We again proceed via induction on $m$. Let $m=1$ and
$r^{\kap(p-1)}>p^{-1}$.

As in the preceding proof put $b_i':=h_i^p-1$ for all $i$. The
elements $\bfh^\beta:=h_1^{\beta_1}\cdot\cdot\cdot h_d^{\beta_d}$
with $\beta\in\{\al\in\N_0^d,~\al_i<p~\forall i \}=:\N_{0,<p}^d$
constitute a system $\Rep'$ of representatives for the cosets in
$G/G^1$. Furthermore, the family of pairwise different nonzero
elements
\[\T:=\{\db'^\al\,\db^\beta\}_{(\al,\beta)\in\N_0^d\times\N_{0,<p}^d}\]
 is an orthogonal basis for $(R_r,||.||_r)$. The latter follows from Lem. \ref{Frommer2} in a way similar
to the preceding proof observing that here, $\iota$ is given by
the bijection $\db'^\al\,\db^\beta\mapsto p\al+\beta$. Since the
elements $\db'^\al,~\al\in\N_0^d$ constitute a topological
$K$-basis for $\Rein_r$ the definition of an orthogonal basis
implies that $R_r$ is a finite and free left $\Rein_r$-module on
the finite basis $\db^\beta,~\beta\in\N_{0,<p}^d.$ Mapping
$\db^\beta\mapsto\bfh^\beta$ induces a left $\Rein_r$-module
isomorphism of $R_r$ onto itself. The latter is thus finite free
over $\Rein_r$ on the basis $\Rep'$ and hence on any basis
consisting of coset representatives for $G/G^1$.

Now let $m>1$ and assume that the result holds true for all
numbers strictly smaller than $m$. By the induction hypothesis we
can assume that \[R_r=\oplus_{h\in\Rep'}\,\Rmin_r\,h\] where
$\Rep'$ is a system of representatives for the cosets in
$G/G^{m-1}$. But $G^m$ appears also as first step in the lower
$p$-series of $G^{m-1}$ and so the induction hypothesis applies to
these two groups as well: if the index $p^{-1}<s<1$ in $p^\Q$
satisfies $s^{\kap p}>s^\kap p^{-1}$ then one has
\[R^{(m-1)}_{[s]}=\oplus_{h\in\Rep''}\Rm_{[s]}h\]
where $\Rep''$ is a system of representatives for
$G^{m-1}/G^{m},~R^{(m-1)}_{[s]}$ denotes the completion of $\Rmin$
via the norm $||.||_s^{(m-1)}$ and $\Rm_{[s]}$ is the closure of
$\Rm$ inside this completion.

Now choose $s=r^{p^{m-1}}$. Then we obtain on the one hand
$s^{\kap p}=r^{\kap p^m}>r^\kap p^{-1}>s^\kap p^{-1}$. On the
other hand Prop. \ref{restriction} implies that $||.||_r$
restricts on $\Rmin$ to the norm $||.||^{(m-1)}_{r^{p^{m-1}}}$.
Hence \[
R^{(m-1)}_{[r^{p^{m-1}}]}=R^{(m-1)}_r,~~~\Rm_{[r^{p^{m-1}}]}=\Rm_r\]
and the proof is complete.
\end{proof}
Now assume that $G$ is a locally $L$-analytic group such that its
scalar restriction $G_0$ is uniform. Endowing each $G^m$ with the
locally $L$-analytic structure as an open subgroup of $G$ we
obtain, again by functoriality, a series $(D(G^m,K))_m$ of
Fr\'echet-Stein algebras the transition maps being topological
embeddings. Write $\nr^{(m)}$ for the set of quotient norms on
$\Dgm$. Again we abbreviate $\nr:=\nr^{(0)}$. Moreover, let us
write $res(.)$ resp. $q(.)$ for restricting a norm (from $\Dgo$ or
$\Dg$ to a subalgebra) resp. forming the quotient norm (with
respect to any of the maps $\Dgom\rightarrow\Dgm$).
\begin{corollary}\label{endlich}
Assume $r^{\kap(p^m-1)}>p^{-1}$ and let $||.||_r$ be given on
$\Dgo$. The two norms $res (\nr)$ and $q\circ res (||.||_r)$
induce the same topology on $\Dgm$. Writing $D_{(r)}(G^m,K)$ for
the completion $\Dgr$ becomes a finite and free (left or right)
module over $D_{(r)}(G^m,K)$ with a basis any set of coset
representatives for $G/G^m$.
\end{corollary}
\begin{proof}
Let $\Rep$ be a system of coset representatives for $G/G^m$. By
the preceding lemma $\Dgor$ is a finite free Banach module over
the noetherian Banach algebra $D_{(r)}(G_0^m,K)$. By \cite{BGR},
3.7.2/3 the $||.||_r$-topology on $\Dgo$ equals the direct sum
topology with respect to $\Dgo=\oplus_{h\in\Rep}\Dgom h$ (where
$\Dgom$ carries the $res (||.||_r)$-topology). Furthermore, one
has $I(G_0,K)=\oplus_{h\in\Rep}I(G_0^m,K)h$ as a direct
consequence of the definitions. Passing to completions and then to
quotients yields the assertions.
\end{proof}

\section{Regularity}

In this section we prove the main result. We start by recalling
the notion of an Auslander regular ring (cf. \cite{LVO}, chap.
III.)

~\\Let $R$ be an arbitrary associative unital ring. For any (left
or right) $R$-module $N$ the {\it grade} $j_R(N)$ is defined to be
either the smallest integer $l$ such that Ext$_R^l(N,R)\neq 0$ or
$\infty$. A left and right noetherian regular ring $R$ is called
{\it Auslander regular} if every finitely generated left or right
$R$-module $N$ satisfies the Auslander condition (AC): for any
$l\geq 0$ and any $R$-submodule $L\subseteq$\,Ext$_R^l(N,R)$ one
has $j_R(L)\geq l$.

A commutative noetherian regular ring is Auslander regular
(\cite{LVO}, III.2.4.3). A complete filtered ring whose graded
ring is Auslander regular of global dimension $d$ is Auslander
regular and has global dimension $\leq d$ (\cite{LVO}, II.2.2.1,
II.3.1.4, III.2.2.5). We deduce
\begin{proposition}\label{subsheaf}
Let $G$ be a locally $L$-analytic group of dimension  $d$ that is
uniform and satifies {\rm (L)}. If $r^\kap<\radius$ then $\Dgr$ is
an Auslander regular ring of global dimension $\leq d$.
\end{proposition}
\begin{proof}
According to Prop. \ref{quotient1} the graded ring $\gr\Dgr$
equals a polynomial ring over $\gor K$ in $d$ variables. Since
$\gor K$ equals Laurent polynomials in one variable we obtain, by
the preceding remarks, the result with the bound $d+1$ on the
global dimension. Using Prop. \ref{result1b} we may deduce, by a
computation similar to \cite{ST5}, (proof of) Lem. 4.8 that $\gr
F^0_r\Dgr$ is isomorphic to a polynomial ring $k[X_0,...,X_d]$
over the residue field $k$ of $K$. Hence, the global dimension of
$F^0_r\Dgr$ is $\leq d+1$ by the above remarks. On the other hand,
multiplication by $p$ is zero on any simple module over
$F^0_r\Dgr$ (\cite{LVO}, I.3.5.5) whence
$\Dgr=\Q_p\otimes_{\Z_p}F^0_r\Dgr$ has global dimension $\leq d$
(\cite{MCR}, 7.4.3/4).
\end{proof}
We recall the following facts from general ring theory.
\begin{lemma}\label{AR1}
Let $R_0\subseteq R_1$ be an extension of unital rings. Suppose
there are units $b_1=1,b_2,...,b_t\in R_1^\times$ which form a
basis of $R_1$ as (left and right) $R_0$-module and which satisfy:

1. $b_iR_0=R_0b_i$ for any $1\leq i\leq t$,

2. for any $1\leq i,j\leq t$ there is $1\leq k\leq t$ such that
$b_ib_j\in b_kR_0$,

3. for any $1\leq i\leq t$ there is $1\leq l\leq t$ such that
$b_i^{-1}\in b_lR_0$.

Suppose $t$ is invertible in $R_0$.

Then: $R_0$ is noetherian if and only if $R_1$ is noetherian. In
this case both rings have the same global dimension.
\end{lemma}
\begin{proof}
E.g. \cite{ST5}, Lem. 8.8.
\end{proof}
\begin{corollary}\label{AR2} Keeping the assumptions $R_0$ is an
Auslander regular ring if and only if this holds true for $R_1$.
\end{corollary}
\begin{proof}
Assume that $R_0$ is Auslander regular. Let $N$ be a finitely
generated left or right $R_1$-module. We have a group isomorphism
\begin{equation}\label{ecu}{\rm~Ext}_{R_1}^\ast(N,R_1)\stackrel{\sim}{\longrightarrow}
{\rm Ext}^\ast_{R_0}(N,R_0).\end{equation} Indeed, using a
projective resolution of $N$ by finitely generated free
$R_1$-modules this reduces to $\ast=0$ and $N=R_1$. Denoting by
$l: R_1\rightarrow R_0$ the projection onto the $b_1$-component
the map $\Phi\mapsto l\circ\Phi$ is the desired bijection. Now let
$L\subseteq$~Ext$_{R_1}^l(N,R_1)$ be any $R_1$-submodule. Since
$R_1$ is noetherian and $N$ is finitely generated so is
Ext$_{R_1}^l(N,R_1)$ and $L$. Consider $L$ as $R_0$-module. Then
it is finitely generated and so from $L\subseteq
$~Ext$_{R_0}^l(N,R_0)$ we deduce by (AC) for the finitely
generated $R_0$-module $N$ that $j_{R_0}(L)\geq l$. But this
implies $j_{R_1}(L)\geq l$ by (\ref{ecu}).

Conversely, suppose that $R_1$ is Auslander regular. Let $N$ be a
finitely generated left $R_0$-module and let $L\subseteq\,
$Ext$_{R_0}^l (N,R_0)$ be any right $R_0$-module. It is finitely
generated by the same argument as above. Put
$N_1:=R_1\otimes_{R_0}N$ whence
\[L\otimes_{R_0}R_1\subseteq
{\rm Ext}_{R_0}^l(N,R_0)\otimes_{R_0}R_1={\rm
Ext}_{R_1}^l(N_1,R_1).\] By (AC) for $N_1$ we have
$j_{R_1}(L\otimes_{R_0}R_1)\geq l$ and so
Ext$_{R_1}^k(L\otimes_{R_0}R_1,R_1)=0$ for all $k<l$. One has
\[{\rm Ext}_{R_1}^k(L\otimes_{R_0}R_1,R_1)=R_1\otimes_{R_0}{\rm
Ext}_{R_0}^k(L,R_0)\] as left $R_1$-modules and so ${\rm
Ext}_{R_0}^k(L,R_0)=0$ for all $k<l$ by faithful flatness of
$R_1$. By definition we obtain $j_{R_0}(L)\geq l$. The proof for
right modules being the same this shows that $R_0$ is Auslander
regular.
\end{proof}
Consider a number $p^{-1}<s<1$ in $p^\Q$ such that
$s^\kap<\radius$. Define $S$ to be the set of all positive real
$p$-power roots of such numbers. Clearly, $S\subseteq
(p^{-1},1)\cap p^\Q$ and, most importantly, $S$ contains a
sequence $s\uparrow 1$.

Turning back to our compact locally $L$-analytic group $G$ which
is uniform satisfying (L) recall that we have a filtration of
$\Dg$ by Fr\'echet-Stein subalgebras
\[\Dg\supseteq D(G^1,K)\subseteq D(G^2,K)\supseteq...\]
coming via functoriality from the lower $p$-series $(G^m)_m$ of
$G$.

It is clear from the definition that the subfamily of norms
$||.||^{(m)}_{\bar{r}},~r\in S$ gives each $\Dgm$ a
Fr\'echet-Stein structure.
\begin{lemma}\label{interim}
Let $r\in S$ and choose $p^{-1}<s<1,~s^\kap<\radius$ such that
$r^{p^{m}}=s$ for some $m\in\N_0$. The norm $res (\nr)$ induces on
$D(G^{m},K)$ the same topology as the norm
$||.||_{\bar{s}}^{(m)}$. Furthermore, $\Dgr$ is a finite and free
(left or right) module over the subring $D_s(G^{m},K)$ on a basis
any set of coset representatives for $G/G^{m}$. In particular,
$\Dgr$ is Auslander regular if and only if this is true for
$D_s(G^{m},K)$ and both rings have the same global dimension.
\end{lemma}
\begin{proof}
The last two claims follow from the second according to Lem.
\ref{AR1} and Cor. \ref{AR2}.

Let us  assume $p\neq 2$ so that $\kap=1$. Now $p^{-1}<s=r^{p^m}$
implies $p^{-1}<r^{(p^m-1)}$. Thus, Cor. \ref{endlich} yields that
$res (\nr)$ induces the same topology as $q\circ res (||.||_r)$ on
the subring $\Dgm\subseteq\Dg$. Moreover, $\Dgr$ is a finite free
module over the closure $D_{(r)}(G^m,K)$ with basis any set of
coset representatives for $G/G^m$. But, according to
Prop.~\ref{restriction} $res (||.||_r)=||.||_{r^{p^m}}^{(m)}$ on
$\Dgom$. Thus, $res (\nr)$ induces in fact the same topology as
$q(||.||_{r^{p^m}}^{(m)})=||.||_{\bar{s}}^{(m)}$ on $D(G^{m},K)$
whence $D_{(r)}(G^m,K)=D_s(G^{m},K)$.

This settles the case $p\neq 2$ and the case $p=2$ follows in a
similar fashion.
\end{proof}
\begin{proposition}
\label{result2} Consider the usual family $\nr,~p^{-1}<r<1,~r\in
p^\Q$ of quotient norms on $\Dg$. If $r\in S$ then $\Dgr$ is
Auslander regular and of global dimension $\leq d$.
\end{proposition}
\begin{proof}
Let $r\in S$. Choose $m\in\N_0$ such that $r^{p^{m}}=:s$ satisfies
$p^{-1}<s<1$ in $p^\Q$ with $s^\kap<\radius$. Then $D_s(G^m,K)$ is
an Auslander regular ring of global dimension $\leq d$ according
to Prop. \ref{subsheaf}. Thus, by the preceding lemma the same is
true for $\Dgr$.
\end{proof}
\begin{theorem}
\label{result3}Let $G$ be a compact locally $L$-analytic group of
dimension $d$. Then $\Dg$ admits a two-sided $K$-Fr\'echet-Stein
structure consisting of Auslander regular Banach algebras of
global dimension $\leq d$.
\end{theorem}
\begin{proof}
Choose an open normal subgroup $H$ which is uniform and satisfies
(L). Consider the set of norms $||.||_r,~r\in S$ on $\Do$ and let
$q(||.||_r),~r\in S$ be the family of quotient norms on $\Dh$.
Choose a system $\Rep$ of representatives $g_i$ containing $1$ for
the cosets in $G/H$ and use the decomposition $\Dgo=\oplus_i\,\Do
g_i$ to define on $\Dgo$ for each $r\in S$ the maximum norm:
\[||\sum_i\lambda_i\,g_i||_r:=\max_i\,||\lambda_i||_r.\]
As explained in section 3 the corresponding completions $\Dgor$
resp. their quotients $\Dgr$ constitute a two-sided
$K$-Fr\'echet-Stein structure on $\Dgo$ resp. $\Dg$. Now a direct
calculation gives that such a quotient norm $\nr$ on $\Dg$
restricts to $q(||.||_r)$ on the subring $\Dh$. Hence, the ring
extension $\Dhr\subseteq\Dgr$ satisfies the assumptions of
Cor.~\ref{AR2}. But the preceding proposition yields that $\Dhr$
is Auslander regular of global dimension $\leq d$.
\end{proof}

\section{Dimension and duality}

In this last section $G$ denotes an {\it arbitrary} locally
$L$-analytic group. We indicate two applications of our results.

For any compact open subgroup $H\subseteq G$ the algebra $\Dh$ is
endowed with the two-sided Fr\'echet-Stein structure exhibited in
the last theorem. It then satisfies the axiom {\rm (DIM)}
formulated in \cite{ST5}, Sect. 8. In [loc.cit.], Sect. 6 the
authors introduce the abelian category $\C_G$ of {\it
coadmissible} (left) $\Dg$-modules. We do not recall the precise
definition here. As an immediate consequence of {\rm (DIM)} the
codimension theory, as developed in [loc.cit.] over $\Q_p$, is now
available on $\C_G$. In particular, any nonzero $M\in\C_G$ has a
well-defined codimension $j_{D(H,K)}(M)$, independent of the
choice of $H$ and bounded above by dim$_L\,G$. A coadmissible
module which is either zero or has maximal codimension equal to
dim$_L G$ is called {\it zero-dimensional}.

As an application we prove that coadmissible modules coming from
smooth or, more general, $U(\lie)$-finite $G$-representations (as
studied in \cite{ST1}) are zero-dimensional. These are natural
generalizations of the theorems \cite{ST5}, 8.12 and 8.15. For the
basic properties of $\C_G$ and its codimension which will be used
in the following we refer to [loc.cit.].
\begin{proposition}
Consider a $d$-dimensional locally $L$-analytic group $H$ which is
uniform and satisfies {\rm (L)} with corresponding $L$-basis
$\x_1,...,\x_d$ of $\lie_L$, the Lie algebra of $H$. Let
$\lambda_1,...,\lambda_d$ be elements of $U(\lie_L)$, the
enveloping algebra, such that $\lambda_j=P_j(\x_j)$ where $P_j$ is
a nonzero polynomial in $L[X]$ with $P_j(0)=0$. Denote by $J$ the
left ideal of $\Dh$ generated by the $\lambda_j$. The coadmissible
module $\Dh/J$ is zero-dimensional.
\end{proposition}
\begin{proof}
Since $J$ is finitely generated the $\Dh$-modules $J$ and $\Dh/J$
are coadmissible. It thus suffices to fix $r\in S$ and prove
\[j_{\Dhr}\,(\Dhr/\Dhr J)\geq d.\] By Lem. \ref{interim} we have with
$R_1:=\Dhr,~R_0:=D_s(H^m,K)$ that $R_1$ is free as $R_0$-module on
the finite basis $\mathcal{R}$. Here, $m\in\N_0$ is appropriately
chosen, $s=r^{p^m}$ and $\mathcal{R}$ is a finite system of
representatives for the cosets in $H/H^m$. According to (the proof
of) Cor. \ref{axiomLm} the uniform group $H^m$ satisfies (L) with
the $L$-basis $p^m\x_j$. Let $h_{11},...,h_{nd}$ denote the
induced minimal ordered generating system for $H^m$ and put as
usual $b_j:=h_{1j}-1$. According to Prop. \ref{result1b} the
elements $\db^\al,~\al\in\N_0^d$ are a topological $K$-basis for
the Banach space $R_0$ and moreover, orthogonal with respect to
$||.||^{(m)}_{\bar{s}}$. As graded ring we have $\grd R_0=(\gor
K)[\sigma(b_1),...,\sigma(b_d)]$ where $\sigma$ denotes the
principal symbol map.

 Since $H^m\subseteq H$ is open we may consider the left ideal $J_0$
generated by $\lambda_1,...,\lambda_d$ inside $R_0$. Abbreviate
$J_1:=\Dhr J$ and use faithful flatness of $R_1$ over $R_0$ to
obtain $j_{R_1}(R_1/J_1)=j_{R_0}(R_0/J_0).$ Since $R_0$ is Zariski
with Auslander regular graded ring the right-hand side equals
$j_{\grd R_0}(\grd R_0/\grd J_0)$ according to [\cite{LVO},
III.2.5.2. By standard commutative algebra the commutative
noetherian regular and catenary domain $\grd R_0$ is
Cohen-Macaulay with respect to Krull dimension (on the category of
finitely generated modules) (\cite{BH}, Cor. 3.5.11). Since $\grd
R_0$ has Krull dimension $d+1$ it therefore suffices to show that
the finitely generated module $\grd R_0/\grd J_0$ has Krull
dimension $\leq 1$.

Now identifying the Lie algebras of $H^m$ and $H_0^m$ over $\Q_p$
we have $p^m\x_j=\log (1+b_j)$ in $R_0$. In $\grd R_0$ we obtain
$\sigma(\log(1+b_j))=\sigma(b_j)$ by definition of the norm
$||.||^{(m)}_{\bar{s}}$ and the fact that $s^\kap<\radius$. Since
$\sigma(p^m)\in\gor L$ is a unit and since $P_j(0)=0$ we obtain
\[\sigma(\lambda_j)=\sigma(P_j(\x_j))=P_j'(\sigma(b_j))\]
with some nonconstant polynomial $P_j'\in (\gor L)[X]$. Since
$\sigma(\lambda_j)\in\grd J_0$ we have a surjection
\[(\gor K)[\sigma(b_1)]/(\sigma(\lambda_1))\otimes_{\gor K}...\otimes_{\gor K}(\gor
K)[\sigma(b_d)]/(\sigma(\lambda_d))\rightarrow\grd R_0/\grd J_0.\]
Each $(\gor K)[\sigma(b_j)]/(\sigma(\lambda_j))$ is finitely
generated as a module over the one-dimensional ring $\gor K$ and
hence,~so is $\grd R_0/\grd J_0$. Then $\grd R_0/\grd J_0$ has
Krull dimension $\leq 1$.
\end{proof}
\begin{theorem}
Let $M\in\C_G$. If the action of the universal enveloping algebra
$U(\lie_L)$ is locally finite, i.e., if $U(\lie_L)x,$ for any
$x\in M$, is a finite dimensional $L$-vectorspace then $M$ is
zero-dimensional.
\end{theorem}
\begin{proof}
Fix an open subgroup $H\subseteq G$ that is uniform satisfying (L)
and let $x\in M$. It suffices to check that the cyclic module $\Dh
x\in\C_H$ is zero-dimensional. Let $\x_1,...,\x_d$ be a basis of
$\lie_L$ realising condition (L) for $H$. Write
$D(H,K)x=D(H,K)/J_1$ with some left ideal $J_1$ of $D(H,K)$ and
let $J_0$ be the kernel ideal of $U(\lie_L)\rightarrow
U(\lie_L)x$. Then $J_0\subseteq J_1$ and the left ideal $J_0$ of
$U(\lie_L)$ has, by assumption, finite codimension in $U(\lie_L)$.
Hence, there is a left ideal $J\subseteq J_1$ in $D(H,K)$
satisfying the assumptions of the preceding proposition. Thus
$D(H,K)/J$ is zero-dimensional and hence, so is its quotient
$D(H,K)/J_1$.
\end{proof}
A smooth $G$-representation $V$ is called {\it admissible-smooth}
if, for any compact open subgroup $H\subseteq G$, the vector
subspace $V^H$ of $H$-invariant vectors in $V$ is finite
dimensional. An admissible-smooth $G$-representation equipped with
the finest locally convex topology is admissible in the sense that
its strong dual lies in $\C_G$ (\cite{ST5}, Thm. 6.6).
\begin{corollary}\label{smoothzero}
If $V$ is an admissible-smooth $G$-representation then the
corresponding coadmissible $D(G,K)$-module is zero-dimensional.
\end{corollary}
\begin{proof}
By smoothness the derived Lie algebra action on $V$ is trivial.
The corresponding coadmissible module $M$ is the strong dual
$V'_b$ with $D(G,K)$-action induced by the contragredient
$G$-action. Thus $U(\lie_L)x=Lx$ for all $x\in M$ and the
preceding theorem applies.
\end{proof}
\begin{theorem}
Let $G$ be compact and $M\in\C_G$. If the $K$-vectorspace
\[M_r:=\Dgr\otimes_{\Dg} M\] is finite dimensional for all $r\in
S$ then $M$ is zero-dimensional.
\end{theorem}
\begin{proof}
Let $H\subseteq G$ be an open uniform subgroup satisfying (L).
Letting $x\in M$ it suffices to prove that $D(H,K)x\in\C_H$ is
zero-dimensional. Write $\Dh x=D(H,K)/J,~J$ some left ideal and
fix $r\in S$. Since $\Dhr$ is flat over $\Dh$ according to
\cite{ST5}, Remark 3.2 we have an inclusion of $\Dhr$-modules
\[\Dhr/\Dhr J=\Dhr\otimes_{\Dh} \Dh x\longrightarrow M_r\] and hence the left-hand
side is finite-dimensional over $K$. If $\lie_L$ denotes the Lie
algebra the left ideal $J_0:=(K\otimes_L U(\lie_L))\cap (\Dhr J)$
has therefore finite $K$-codimension in $K\otimes_LU(\lie_L)$. As
in the proof of the preceding theorem $\Dh/\Dh J_0$ and also
$\Dh/J$ are zero-dimensional.
\end{proof}

As a second consequence let us remark that the duality theory for
admissible locally analytic representations (in the sense of
\cite{ST6}) is now completely available over the base field $L$.
Thus, also for locally $L$-analytic groups $G$, the duality
functor (defined on the bounded derived category of $\Dg$-modules
with coadmissible cohomology) is an anti-involution. Due to the
presence of the codimension the category $\mathcal{C}_G$ is
filtered by abelian subquotient categories and the functor is
computed as a particular Ext-group on each of them.

{\it Acknowledgements.}~This work is based on the author's thesis.
He would like to thank his thesis advisor Peter Schneider for
guidance and support. He is also grateful to Matthias Strauch and
Jan Kohlhaase for many helpful discussions.

\end{document}